\documentclass[12pt]{amsart}

\theoremstyle{plain}    
\newtheorem{thm}{Theorem} 
\numberwithin{figure}{section} 
\theoremstyle{plain}    
\newtheorem{lemma}[thm]{Lemma} 
\newtheorem{prop}[thm]{Proposition}

\newtheorem{defi}[thm]{Definition}

\usepackage{amscd,amssymb,comment,euscript}
\newcount\theTime
\newcount\theHour
\newcount\theMinute
\newcount\theMinuteTens
\newcount\theScratch
\theTime=\number\time
\theHour=\theTime
\divide\theHour by 60
\theScratch=\theHour
\multiply\theScratch by 60
\theMinute=\theTime
\advance\theMinute by -\theScratch
\theMinuteTens=\theMinute
\divide\theMinuteTens by 10
\theScratch=\theMinuteTens
\multiply\theScratch by 10
\advance\theMinute by -\theScratch

\def\today{{\number\day\space
 \ifcase\month\or
  January\or February\or March\or April\or May\or June\or
  July\or August\or September\or October\or November\or December\fi
 \space\number\year}}

\newcommand\Ac{{\mathcal{A}}}

\newcommand\Cpx{{\mathbf C}}

\newcommand\Fb{{\mathbf F}}

\newcommand\fdim{\text{\rm fdim}\,}

\newcommand\Mc{{\mathcal{M}}}

\newcommand\Mct{{\widetilde\Mc}}

\newcommand\Nats{{\mathbf N}}

\newcommand\Nc{{\mathcal{N}}}

\newcommand\nt{{\tilde n}}

\newcommand\Qc{{\mathcal{Q}}}

\newcommand\Qct{{\widetilde\Qc}}

\newcommand\rt{{\tilde r}}

\begin{document}

\pagestyle{myheadings}

\title{Rescalings of free products of II$_1$--factors}
 
\author{Ken Dykema}
\address{\hskip-\parindent
Department of Mathematics\\
Texas A\&M University\\
College Station TX 77843--3368, USA}
\email{Ken.Dykema@math.tamu.edu}

\author{Florin R\u adulescu}

\address{\hskip-\parindent
Florin R\u adulescu \\
Department of Mathematics \\
University of Iowa \\
Iowa City IA 52242--1466, USA}
\email{radulesc@math.uiowa.edu}

\thanks{K.D.\ supported in part by NSF grant DMS--0070558.
F.R.\ supported in part by NSF grantDMS--9970486.
Both authors thank also the Mathematical Sciences Research Institute.
Research at MSRI is supported in part by NSF grant DMS--9701755.}

\date{19 March 2001}

\begin{abstract}
We introduce notation $\Qc(1)*\cdots*\Qc(n)*L(\Fb_r)$ for von Neumann algebra II$_1$--factors
where $r$ is allowed to be negative.
This notation is defined by rescalings of
free products of II$_1$--factors, and is proved to be consistent with
known results and natural operations.
We also give two statements which we prove are equivalent to isomorphism
of free group factors.
\end{abstract}

\maketitle

\markboth{\tiny Rescalings of free products}{\tiny Rescalings of free products}
 
\section*{Introduction}
The rescaling $\Mc_t$ of a II$_1$--factor $\Mc$ by a positive number $t$
was introduced by Murray and von Neumann~\cite{MvN}.
In the paper~\cite{DR}, we showed that if $\Qc(1),\ldots,\Qc(n)$
are II$_1$--factors ($n\in\{2,3,\ldots\}$) and if $0<t<\sqrt{1-1/n}$ then
\begin{equation}
\label{eq:DR}
\big(\Qc(1)*\cdots*\Qc(n)\big)_t
\cong\Qc(1)_t*\cdots*\Qc(n)_t*L(\Fb_r),
\end{equation}
where $r=(n-1)(t^{-2}-1)$.
Here $L(\Fb_r)$, $r>1$, is an interpolated free group factor~(\cite{D94}, \cite{R}).
For $\sqrt{1-1/n}\le t<1$, we proved a similar formula, where $L(\Fb_r)$
was replaced by a hyperfinite von Neumann algebra with specified tracial state
having free dimension~(\cite{D93}) equal to $r=(n-1)(t^{-2}-1)\le1$.
If one tries to use the formula~\eqref{eq:DR}
when $t>1$, one obtains $L(\Fb_r)$ with $r<0$.

In this note we introduce notation
\begin{equation}
\label{eq:Notation}
\Qc(1)*\cdots*\Qc(n)*L(\Fb_r)\qquad (n\in\Nats,\,1-n<r\le\infty).
\end{equation}
If $r>1$ then $L(\Fb_r)$ in~\eqref{eq:Notation} is an interpolated
free group factor, while if $r\le1$ then~\eqref{eq:Notation}
defines a II$_1$ factor which is the rescaling by $t$ of $\Qc(1)_{1/t}*\cdots*\Qc(n)_{1/t}$
if $n=2$ or of $\Qc(1)_{1/t}*L(\Fb_2)$ if $n=1$
for an appropriate $t>1$.
We will prove that this notation is consistent with known results and natural operations
involving free products.
The notation~\eqref{eq:Notation}
provides an elegent means of describing rescalings of free products
of II$_1$--factors,
and is used in~\cite{DH} and~\cite{D}.

Finally, we show that if the free group factors are isomorphic to
each other then $\Qc(1)*\Qc(2)\cong\Qc(1)*\Qc(2)*L(\Fb_\infty)$
for all II$_1$--factors $\Qc(1)$ and $\Qc(2)$ and we give
one additional equivalent condition.
It is conceivable that these conditions may be used to prove
nonisomorphism of free group factors.

\section*{Rescalings}
\vskip1ex
Recall that the interpolated free group factors rescale as follows:
\begin{equation}
\label{eq:LFres}
L(\Fb_r)_t\cong L(\Fb_{1+t^{-2}(r-1)}),\qquad(1<r\le\infty,\,0<t<\infty),
\end{equation}
(see~\cite{D94}, \cite{R}).

\begin{lemma}
\label{lem:1}
Let $n\in\Nats$, let $\Qc(1),\ldots,\Qc(n)$ be II$_1$--factors, let $1<r\le\infty$
and let
\[
\Mc=\Qc(1)*\cdots*\Qc(n)*L(\Fb_r).
\]
Then for every $0<t<\sqrt{1+(r-1)/n}$,
\[
\Mc_t\cong\Qc(1)_t*\cdots*\Qc(n)_t*L(\Fb_{t^{-2}r+(n-1)(t^{-2}-1)}).
\]
\end{lemma}
\begin{proof}
If $t\le1$ then this follows from~\cite{DR}; see~\eqref{eq:DR} and~\eqref{eq:LFres} above.
Suppose $t>1$.
Note that $t$ is taken so that $t^{-2}r+(n-1)(t^{-2}-1)>1$.
Applying~\eqref{eq:DR} and~\eqref{eq:LFres}, we have
\[
\Big(\Qc(1)_t*\cdots*\Qc(n)_t*L(\Fb_{t^{-2}r+(n-1)(t^{-2}-1)})\Big)_{\frac1t}
\cong\Qc(1)*\cdots*\Qc(n)*L(\Fb_r).
\]
\end{proof}

\begin{prop}
\label{prop:*LFr}
Let $n\in\Nats$, let $\Qc(1),\ldots,\Qc(n)$ be II$_1$--factors
and let
\[
1-n<r\le1.
\]
Then there is a II$_1$--factor $\Mc$, unique up to isomorphism, such that
\[
\Mc_t\cong\Qc(1)_t*\cdots*\Qc(n)_t*L(\Fb_{t^{-2}r+(n-1)(t^{-2}-1)})
\]
whenever $0<t<\sqrt{1+(r-1)/n}$.
\end{prop}
\begin{proof}
Let $0<s<t<\sqrt{1+(r-1)/n}$ and let $\Mc$ and $\Mct$ be II$_1$--factors
such that
\begin{align*}
\Mc_s&\cong\Qc(1)_s*\cdots*\Qc(n)_s*L(\Fb_{s^{-2}r+(n-1)(s^{-2}-1)}) \\
\Mct_t&\cong\Qc(1)_t*\cdots*\Qc(n)_t*L(\Fb_{t^{-2}r+(n-1)(t^{-2}-1)}).
\end{align*}
Then using Lemma~\ref{lem:1} we have
\[
\Mct_s=(\Mct_t)_{\frac st}
\cong\Qc(1)_s*\cdots*\Qc(n)_s*L(\Fb_{s^{-2}r+(n-1)(s^{-2}-1)})\cong\Mc_{s.}
\]
\end{proof}

\begin{defi}\rm
\label{def:*LFr}
We denote the unique factor $\Mc$ in Proposition~\ref{prop:*LFr} by
\begin{equation*}
\Qc(1)*\cdots*\Qc(n)*L(\Fb_r).
\end{equation*}
\end{defi}

\begin{prop}
\label{prop:properties}
Let $n\in\Nats$, let $\Qc(1),\Qc(2),\ldots,\Qc(n)$ be II$_1$--factors
and let $1-n<r\le\infty$.
\renewcommand{\labelenumi}{(\roman{enumi})}
\begin{enumerate}

\item
If $0<t<\infty$ then
\[
\big(\Qc(1)*\cdots*\Qc(n)*L(\Fb_r)\big)_t
\cong\Qc(1)_t*\cdots*\Qc(n)_t*L(\Fb_{t^{-2}r+(n-1)(t^{-2}-1)}).
\]

\item If $\sigma$ is a permutation of $\{1,2,\ldots,n\}$ then
\[
\Qc(1)*\cdots*\Qc(n)*L(\Fb_r)\cong\Qc(\sigma(1))*\cdots*\Qc(\sigma(n))*L(\Fb_r).
\]

\item If $\Qc(1)=L(\Fb_s)$ with $1<s\le\infty$ then
\[
\Qc(1)*\cdots\Qc(n)*L(\Fb_r)\cong\begin{cases}
\Qc(2)*\cdots*\Qc(n)*L(\Fb_{r+s})&\text{if }n\ge2\\
L(\Fb_{r+s})&\text{if }n=1.
\end{cases}
\]

\item If $n\ge2$ and if $r>2-n$ then
\[
\Qc(1)*\cdots*\Qc(n)*L(\Fb_r)\cong\Qc(1)*\big(\Qc(2)*\cdots*\Qc(n)*L(\Fb_r)\big)_.
\]

\item If $\Qc(1)=\Nc(1)*\Nc(2)$ where $\Nc(1)$ and $\Nc(2)$ are II$_1$--factors then
\[
\Qc(1)*\cdots*\Qc(n)*L(\Fb_r)\cong\Nc(1)*\Nc(2)*\Qc(2)*\cdots*\Qc(n)*L(\Fb_r).
\]

\item If $0<\rt\le\infty$ then
\[
\big(\Qc(1)*\cdots*\Qc(n)*L(\Fb_r)\big)*L(\Fb_\rt)
\cong\Qc(1)*\cdots*\Qc(n)*L(\Fb_{r+\rt}).
\]

\item If $\nt\in\Nats$, if $\Qct(1),\ldots,\Qct(\nt)$ are II$_1$--factors
and if $1-\nt<\rt\le\infty$ then
\begin{gather*}
\big(\Qc(1)*\cdots*\Qc(n)*L(\Fb_r)\big)*
\big(\Qct(1)*\cdots*\Qct(\nt)*L(\Fb_\rt)\big)\cong \\
\cong\Qc(1)*\cdots*\Qc(n)*\Qct(1)*\cdots*\Qct(\nt)*L(\Fb_{r+\rt}).
\end{gather*}

\item If $n\ge2$ then
\[
\Qc(1)*\cdots*\Qc(n)*L(\Fb_0)\cong\Qc(1)*\cdots*\Qc(n).
\]

\item If $\Nc$ is a II$_1$--factor and if $\Ac$ is a von Neumann algebra with
specified normal faithful tracial state, where $\Ac\ne\Cpx$ and $\Ac$ is either finite
dimensional, hyperfinite, an interpolated free group factor or
a (possibly countably infinite) direct sum of these then
\[
\Nc*\Ac\cong\Nc*L(\Fb_r)
\]
where $r=\fdim(\Ac)$ is the free dimension of $\Ac$ (see~\cite{D93}).
\end{enumerate}
\end{prop}
\begin{proof}
For~(i), if $r>1$ then this is Lemma~\ref{lem:1}.
If $r\le1$ but $0<t<\sqrt{1+(r-1)/n}$ then this is Definition~\ref{def:*LFr}.
Suppose $r\le1$ and $\sqrt{1+(r-1)/n}\le t<\infty$.
Let $\lambda>0$ be such that $\lambda t<\sqrt{1+(r-1)/n}$.
Then applying Definition~\ref{def:*LFr} twice gives
\begin{align*}
\Big(\Qc(1)*\cdots*\Qc(n)*L(\Fb_r)\Big)_{\lambda t}
&\cong\Qc(1)_{\lambda t}*\cdots *\Qc(n)_{\lambda t}
 *L(\Fb_{\lambda^{-2}t^{-2}r+(n-1)(\lambda^{-2}t^{-2}-1)}) \\
&\cong\Big(\Qc(1)_t*\cdots*\Qc(n)_t*L(\Fb_{t^{-2}r+(n-1)(t^{-2}-1)})\Big)_{\lambda.}
\end{align*}

Now the proofs of (ii)--(viii) are obtained by rescaling both sides
of the desired isomorphisms by the same $t>0$ which is small enough
and applying~(i) and perhaps equation~\eqref{eq:DR}.
For example, to prove~(vii) let 
\[
0<t<\min\Big(\frac1{\sqrt2},\sqrt{1+\frac{r-1}n},\sqrt{1+\frac{\rt-1}\nt},
\sqrt{1+\frac{r+\rt-1}{n+\nt}}\Big)
\]
and use~(i) three times to get
\begin{align*}
\Big(\big(\Qc(1)*&\cdots*\Qc(n)*L(\Fb_r)\big)
 *\big(\Qct(1)*\cdots*\Qct(\nt)*L(\Fb_\rt)\big)\Big)_t \\
&\cong\Big(\Qc(1)*\cdots*\Qc(n)*L(\Fb_r)\Big)_t
 *\Big(\Qct(1)*\cdots*\Qct(\nt)*L(\Fb_\rt)\Big)_t*L(\Fb_{t^{-2}-1}) \\
&\cong\Qc(1)_t*\cdots*\Qc(n)_t*\Qct(1)_t*\cdots*\Qct(\nt)_t
 *L(\Fb_{t^{-2}(r+\rt)+(n+\nt-1)(t^{-2}-1)}) \\
&\cong\Big(\Qc(1)*\cdots*\Qc(n)*\Qct(1)*\cdots*\Qct(\nt)*L(\Fb_{r+\rt})\Big)_{t.}
\end{align*}

For~(ix), if $k\in\Nats$ is large enough then by~\cite{D93}, $M_k(\Cpx)*\Ac$
is the interpolated free group factor $L(F_{r+1-k^{-2}})$.
By~\cite[Thm. 1.2]{D93},
\begin{align*}
(\Nc*\Ac)_{\frac1k}&\cong\big((\Nc_{\frac1k}\otimes M_k(\Cpx))*\Ac)_{\frac1k}
 \cong\Nc_{\frac1k}*\big(M_k(\Cpx)*\Ac)_{\frac1k}\cong \\
&\cong\Nc_{\frac1k}*L(\Fb_{r+1-k^{-2}})_{\frac1k}\cong\Nc_{\frac1k}*L(\Fb_{k^2r})
 \cong\big(\Nc*L(\Fb_r)\big)_{{\frac1k}_.}
\end{align*}
\end{proof}

Formula~\eqref{eq:DR} can now be extended to all values of $t$.
\begin{thm}
Let $n\in\{2,3,\ldots\}$, let $\Qc(1),\ldots,\Qc(n)$ be II$_1$--factors
and let $0<t<\infty$.
Then
\begin{equation*}
\big(\Qc(1)*\cdots*\Qc(n)\big)_t
=\Qc(1)_t*\cdots*\Qc(n)_t*L(\Fb_{(n-1)(t^{-2}-1)}).
\end{equation*}
\end{thm}
\begin{proof}
Use part~(viii) followed by part~(i) of Proposition~\ref{prop:properties}.
\end{proof}

We know from~\cite{R} (see also~\cite{D94}) that the interpolated free group
factors $(L(\Fb_t))_{1<t\le\infty}$ are either all isomorphic to each other or
all mutually nonisomorphic.

\begin{thm}
The following are equivalent:
\renewcommand{\labelenumi}{(\alph{enumi})}
\begin{enumerate}
\item $L(F_s)\cong L(F_t)$ for some, and then for all, $1<s<t\le\infty$;
\item for every II$_1$--factor $\Qc$ and every $r>0$,
\begin{equation*}
\Qc*L(\Fb_r)\cong\Qc*L(\Fb_\infty);
\end{equation*}
\item for all II$_1$--factors $\Qc(1)$ and $\Qc(2)$,
\begin{equation*}
\Qc(1)*\Qc(2)\cong\Qc(1)*\Qc(2)*L(\Fb_\infty).
\end{equation*}
\end{enumerate}
\end{thm}
\begin{proof}
For (a)$\implies$(b), if $0<t<\sqrt r$ then by part~(i)
of Proposition~\ref{prop:properties},
\[
\big(\Qc*L(\Fb_r)\big)_t\cong\Qc_t*L(\Fb_{t^{-2}r})
\cong\Qc_t*L(\Fb_\infty)\cong\big(\Qc*L(\Fb_\infty)\big)_{t,}
\]
while (b)$\implies$(a) can be seen be choosing $\Qc=L(\Fb_2)$.
For (a)$\implies$(c), if $0<t<1/\sqrt2$ then using Lemma~\ref{lem:1},
\begin{align*}
\big(\Qc(1)*\Qc(2)\big)_t&\cong\Qc(1)_t*\Qc(2)_t*L(\Fb_{t^{-2}-1})\cong \\
&\cong\Qc(1)_t*\Qc(2)_t*L(\Fb_\infty)
\cong\big(\Qc(1)*\Qc(2)*L(\Fb_\infty)\big)_{t.}
\end{align*}
Taking $\Qc(1)\cong\Qc(2)\cong L(\Fb_2)$ shows (c)$\implies$(a).
\end{proof}

\newpage

\bibliographystyle{plain}

\begin{thebibliography}{9}

\bibitem{D93} K.\ Dykema,
{\em Free products of hyperfinite von Neumann algebras and free dimension,}
Duke Math.\ J.\ {\bf 69} (1993), 97-119.

\bibitem{D94} \rule{3em}{.1mm},
{\em Interpolated free group factors,}
Pacific J.\ Math.\ {\bf 163} (1994), 123-135.

\bibitem{D} \rule{3em}{.1mm},
{\em Free subproducts and free scaled products of II$_1$--factors,}
preprint

\bibitem{DH} K.\ Dykema, U.\ Haagerup,
{\em Decomposability of Voiculescu's circular operator and DT--operators,}
in preparation.

\bibitem{DR} K.\ Dykema, F.\ R\u adulescu,
{\em Compressions of free products of von Neumann algebras,}
Math.\ Ann.\ {\bf 316} (2000), 61-82.

\bibitem{MvN} F.J.~Murray and J.~von Neumann,
{\em Rings of operators.~IV,}
Ann.\ of Math.\ {\bf 44} (1943), 716-808.

\bibitem{R} F.\ R\u adulescu,
{\em Random matrices, amalgamated free products and subfactors of the von Neumann algebra of a free group, of noninteger index,}
Invent.\ Math.\ {\bf 115} (1994), 347-389.

\end{thebibliography}

\end{document}